\documentclass[12pt]{article}

\usepackage{amsgen,amsmath,amstext,amsbsy,amsopn,amsfonts,amssymb}
\usepackage[dvips]{graphicx}

\title{On associative algebras with unity and Lie index $1$.}
\author{Vladimir Dergachev}

\def\buzz#1{{\normalfont{\bf{[#1]}}}}

\def\acr{\nonumber \\}

\def\Stab{\textrm{\normalfont Stab}}

\def\Nil{\textrm{\normalfont Nil}}

\def\equivalent{\Leftrightarrow}
\def\imply{\Rightarrow}
\def\implies{\imply}
\def\intersect{\cap}
\def\orthogonal{\bot}
\def\image{\textrm{\normalfont Im}\,}

\def\extended{0}

\newtheorem{theorem}{Theorem}
\newtheorem{definition}{Definition}

\newtheorem{statement}[theorem]{Statement}
\newtheorem{lemma}[theorem]{Lemma}

\newenvironment{proof}{\begin{trivlist}\item[\hskip%
\labelsep{\bf Proof\quad}]}%
{\hfill\qed\end{trivlist}}
\newcommand{\qed}{{\unskip\nobreak\hfil\penalty50\hskip .001pt \hbox{}
          \nobreak\hfil
          \vrule height 1.2ex width 1.1ex depth -.1ex
           \parfillskip=0pt\finalhyphendemerits=0\medbreak}\rm}

\makeatletter
\newcount\@refno
\def\nextref#1#2#3#4{\advance\@refno\@ne
\if@filesw \immediate\write\@auxout
   {\string\bibcite{#1}{\number\@refno}}\fi   
    {{\number\@refno}.\quad{ #2}, {#3}, { #4}.\hfill\\}}

\newcommand{\references}{
\section*{References}
\frenchspacing
\entries\par}

\newcommand{\entries}{
\nextref{DK}
{Vladimir Dergachev and Alexandre Kirillov}
{\em Index of Lie algebras of Seaweed type}
{to appear in {\em Journal of Lie theory}}
\nextref{D}
{Vladimir Dergachev}
{\em Some properties of index of Lie algebras}
{math.RT/0001042}
\nextref
{D1}
{J.Dixmier}
{``Enveloping algebras''}
{American Mathematical Society, Providence (1996) 1-379 }
\nextref{E1}{Elashvili, A. G}
{\em Frobenius Lie algebras}
{Funktsional. Anal. i Prilozhen. {\bf 16} (1982), 94--95}
\nextref{E2}{---}
{\em On the index of a parabolic subalgebra 
of a semisimple Lie algebra}
{Preprint, 1990}
\nextref{GK}
{I.M.Gelfand and A.A.Kirillov}
{\em Sur les corps li\'es aux alg\'ebres enveloppantes des alg\'ebres de Lie}
{ Publications math\'ematiques {\bf 31} (1966) 5-20}

}

\begin{document}
\maketitle
\tableofcontents

\section{Introduction}

This paper continues the study of associative algebras started in \cite{D}. 
We prove several interesting results concerning the roots of characteristic 
polynomial of associative algebra. Further development yeilds a new way
of analysis of associative algebras which we apply to classification of
associative algebras with  unity and Lie index $1$.

This new approach to analysis of associative algebras is remarkable on several
counts. 

First, unlike in classical theory, we do not deal with characteristics
of particular elements in algebra but we analyse their relationship with 
functionals on the algebra.

Secondly, the method bears most resemblance to the famous theory of commutative
$C^*$-algebras, in which particular importance was attached to multiplicative
functionals. These multiplicative functionals can be considered as rank 1
functionals (since the multiplication table of the algebra evaluated on such 
a functional has rank 1). In our method we study the functionals for which this rank is 
maximal possible.

\section{Multiplication table of associative algebra}

\begin{definition}\buzz{Characteristic polynomial of associative algebra}
Let $\mathfrak A$ be an associative algebra, and let $A$ be the multiplication
table of this algebra in some basis $\left\{e_i\right\}$.
$$
\chi(\lambda,\mu,F)=\det(\lambda A+\mu A^T)
$$
here $F$ is an element of a dual space $\mathfrak A'$, and $A$ is considered
as matrix of functions on $\mathfrak A'$.
\end{definition}

Characteristic polynomial $\chi(\lambda,\mu,F)$ is homogeneous in $\lambda$ and
$\mu$ and thus, for each particular $F$, can be decomposed as product of 
factors
$$
\chi(\lambda,\mu,F)=C\lambda^k\mu^k\prod_i\left(\lambda-\alpha_i \mu\right)
$$

\begin{definition}\buzz{$\alpha$-stabilizer of $F$}
Let $F$ be an element of $\mathfrak A'$. We define
$$
\Stab_F(\alpha):=\left\{a\in \mathfrak A : \forall b \in \mathfrak A\imply F(ab)=\alpha F(ba)\right\}=
\ker (\left.A\right|_F-\alpha \left.A^T\right|_F)
$$
where $\left.A\right|_F$ stands for matrix $A$ evaluated in point $F$.
\end{definition}

Obviously $\left(\lambda -\alpha \mu\right)$ divides $\chi$ if and only if $\Stab_F(\alpha)$
is non-trivial. It is also convenient to introduce the following definitions:

\begin{definition}
Let $F$ be an element of $\mathfrak A'$. We define
$$
\Stab_F(\infty):=\left\{a\in \mathfrak A : \forall b \in \mathfrak A\imply F(ba)=0\right\}=
\ker \left.A^T\right|_F
$$

and (this follows from the definition above)

$$
\Stab_F(0):=\left\{a\in \mathfrak A : \forall b \in \mathfrak A\imply F(ab)=0\right\}=
\ker \left.A\right|_F
$$
\end{definition}

\begin{definition}
Let $F$ be an element of $\mathfrak A'$. We define
$$
\Nil_F=\left\{a\in \mathfrak A : \forall b \in \mathfrak A\imply F(ab)=0 \textrm{ and } F(ba)=0\right\}
$$
\end{definition}

\begin{statement} For all distinct $\alpha$ and $\beta$ we have
$$
\Stab_F(\alpha) \cap \Stab_F(\beta)=\Nil_F
$$
\end{statement}

\begin{definition}\buzz{$Q_F$} We define symmetric bilinear form $Q_F$ on $\Stab_F(1)$
 as follows:
$$
Q_F(a,b)=F(ab)
$$
\end{definition}

\begin{statement} The following properties of an associative algebra are 
equivalent:
\begin{itemize}
\item $\Nil_F$ is trivial
\item the characteristic polynomial $\chi(\lambda,\mu,F)$ is non-zero in 
point $F$.
\end{itemize}
\end{statement}

\begin{statement}\label{Q_F} The following properties of an associative algebra are equivalent:
\begin{itemize}
\item $Q_F$ is non-degenerate
\item the maximum power of $\lambda+\mu$ that divides characteristic polynomial
is equal to the index of the algebra
\end{itemize}
\end{statement}

Note also that this non-degeneracy of $Q_F$ implies that characteristic polynomial
of $\mathfrak A$ is non-zero. As an example, all associative algebras with unity and 
index $1$ have these properties. For proof see \cite{D}.

\begin{statement}\label{stab_times} Let $\alpha$ and $\beta$ be two finite numbers. Then
$$
\Stab_F(\alpha)\cdot \Stab_F(\beta) \subset \Stab_F(\alpha\beta)
$$
\end{statement}
\begin{proof}
Let $a \in \Stab_F(\alpha)$, $b \in \Stab_F(\beta)$ and $c \in \mathfrak A$. Then
$$
F((ab)c)=F(a(bc))=\alpha F((bc)a)=\alpha F((b(ca))=\alpha \beta F((ca)b)=\alpha \beta F(c(ab))
$$
\end{proof}

\begin{statement} The following inclusions are true:
\begin{description}
\item[] $\Stab_F(0) \cdot \Stab_F(\infty) \subset \Nil_F$
\item[] $\Stab_F(\alpha) \cdot \Stab_F(\infty) \subset \Stab_F(\infty)$ (when $\alpha\neq 0$)
\end{description}
\end{statement}

\begin{statement} The following equalities are true:
\begin{description}
\item[] $\dim \Stab_F(0)=\dim \Stab_F(\infty)$
\item[] $\dim \Stab_F(\alpha)=\dim \Stab_F(\frac{1}{\alpha})$ 
\end{description}
\end{statement}
\begin{proof}
Indeed
$$
\dim \ker \left.A\right|_F=\dim \ker \left.A^T\right|_F
$$
and
\begin{eqnarray}
\lefteqn{\dim \ker \left(\left.A\right|_F-\alpha\left.A^T\right|_F\right)=}\acr
&&\qquad=\dim\ker\left(\left.A^T\right|_F-\alpha\left.A\right|_F\right)=\acr
&&\qquad=\dim\ker\left(\left.A\right|_F-\frac{1}{\alpha}\left.A^T\right|_F\right)\nonumber
\end{eqnarray}
\end{proof}

\begin{statement} For all $\alpha$ and $\beta$ such that $\alpha \beta \neq 1$ and
$(\alpha,\beta)\neq (\infty,0)$  we have
$$
\Stab_F(\alpha) \cdot \Stab_F(\beta) \subset \ker F
$$
(or for all $a \in \Stab_F(\alpha)$ and $b \in \Stab_F(\beta)$ we have $F(ab)=0$).
\end{statement}

\begin{statement} $\Stab_F(1)$ is commutative (for generic $F$).
\end{statement}
\begin{proof}\ \\
Indeed the Lie algebra derived from $\Stab_F(1)$ is commutative (see \cite{D1}).
\end{proof}

\section{Classification of 2-dimensional associative algebras}
We will now apply the general theory presented above to the problem of classification
of 2-dimensional algebras.

Such an algebra is either commutative (index is $2$) or non-commutative 
(index is $0$). The first case is easy and has been covered in literature.
Let us now concentrate on the second possibility. 

Each algebra of index $0$ must have non-zero characteristic polynomial that
is not divisible by $\lambda+\mu$ (see statement \ref{Q_F}).

\begin{lemma}
Let $\mathfrak A$ be a 2-dimensional associative algebra with index $0$. Let $F$
be a generic element of $\mathfrak A^*$. Then $\Stab_F(\alpha)=0$ for all 
$\alpha\notin \left\{0,1,\infty\right\}$.
\end{lemma}
\begin{proof}
Assume the converse. In view of statement \ref{Q_F} $\Nil_F$ is trivial and thus
$\mathfrak A=\Stab_F(\alpha)\oplus\Stab_F(1/\alpha)$. Because of statement \ref{stab_times}
the multiplication is trivial ($0$) - which contradicts our assumption that
$\mathfrak A$ is non-commutative.
\end{proof}

Thus we must have that $\dim\Stab_F(0)=\dim\Stab_F(\infty)=1$. Let us choose
$x\in\Stab_F(0)$ and $y\in\Stab_F(\infty)$, $x\neq 0$ and $y\neq 0$.

The multiplication table must look like
$$
\begin{array}{cccc}
\times_{\mathfrak A} & \vline & x & y \\
\hline
x &  \vline & \mu x  & 0  \\
y &  \vline & \beta x+\gamma y  & \nu y  \\
\end{array}
$$

Now we will use the fact that $\mathfrak A$ is associative to obtain the following
equations:

\begin{eqnarray}
\label{eq1}
x(yx)=(xy)x & \imply &  \mu\beta x=0\\
\label{eq2}
y(xx)=(yx)x & \imply & \mu\beta x+\mu\gamma y=\mu\beta x+\gamma(\beta x+\gamma y)  \\
\label{eq3}
y(yx)=(yy)x & \imply & \beta(\beta x+\gamma y)+\nu\gamma y=\nu\beta x+\nu\gamma y\\
\label{eq4}
y(xy)=(yx)y & \imply & 0=\gamma\nu y
\end{eqnarray}

Combining first two equations we derive:
$$
\mu(\beta x+\gamma y)=\gamma (\beta x+\gamma y)
$$

And using equations \ref{eq3} and \ref{eq4} we get:
$$
\beta(\beta x+\gamma y)=\nu(\beta x+\gamma y)
$$

Since $\mathfrak A$ is not commutative $\beta x+\gamma y$ must be non-zero. 
Hence $\mu=\gamma$ and $\beta=\nu$. Also using equation \ref{eq1} we get 
$\mu\nu=0$. This leaves us only two choices - one when $\mu$ is non-zero and
one when $\nu$ is non-zero. By multiplying $x$ or $y$ by appropriate constants
we can normalize non-zero parameter ($\nu$ or $\mu$) to $1$.

Thus the only two non-commutative associative algebras are
$$
\begin{array}{cccc}
\times_{\mathfrak A} & \vline & x & y \\
\hline
x &  \vline & x  & 0  \\
y &  \vline & y  & 0  \\
\end{array}
$$
and
$$
\begin{array}{cccc}
\times_{\mathfrak A} & \vline & x & y \\
\hline
x &  \vline & 0  & 0  \\
y &  \vline & x  & y  \\
\end{array}
$$

\section{Classification of 3-dimensional associative algebras with unity}

We will now approach the problem of 
classification of 3-dimensional associative algebras with unity.

For these algebras the index can attain only 2 values: 1 or 3. 

\subsection{Index equals 3}

In this case the algebra is commutative.

\subsection{Index equals 1}

\begin{lemma}
Let $\mathfrak A$ be a 3-dimensional associative algebra with unity and index $1$. Let $F$
be a generic element of $\mathfrak A^*$. Then $\Stab_F(\alpha)=0$ for all 
$\alpha\notin \left\{0,1,\infty\right\}$.
\end{lemma}
\begin{proof}
Indeed, assume that the converse is true. Then there exist non-zero elements $x$ and $y$
such that $x\in \Stab_F(\alpha)$, $y\in \Stab_F(1/\alpha)$, 
$\alpha \notin \left\{0,1,\infty\right\}$ and $F$ is generic. 
Let us normalize $F$ so that $F(1)=1$.

Thus $\mathfrak A$ is spanned by $1$, $x$ and $y$. Since $\alpha^2\neq \alpha$
we must have that $x^2=y^2=0$. Since 
$\Stab_F(\alpha)\cdot \Stab_F(1/\alpha)\subset \Stab_F(1)$ we have 
$xy=\gamma\cdot 1$. 

But 
$$
xyx=(xy)x=(\gamma\cdot 1) \cdot x=x\cdot( \gamma \cdot 1)=x(xy)=x^2y=0
$$

Thus $xy=0$. Analogously $yx=0$. But this contradicts the assumption that our
algebra is not commutative (index equals $1$).
\end{proof}

In view of statement \ref{Q_F} the characteristic polynomial must be non-zero
and divisible only by first power of $\lambda+\mu$. Hence the only other possibility
is for $\dim \Stab_F(0)=\dim \Stab_F(\infty)=1$.

Let $x \in \Stab_F(0)$ and $y\in\Stab_F(\infty)$ such that $x\neq 0$ and $y \neq 0$.
Then $1$, $x$ and $y$ span $\mathfrak A$. Since $\Nil_F=0$ we must have $xy=0$.
Thus the multiplication table should look like:
$$
\begin{array}{ccccc}
\times_{\mathfrak A}&\vline  & 1 &x &y  \\
\hline
1&\vline &1 &x &y  \\
x&\vline &x &\alpha x &0  \\
y&\vline &y &\gamma+\mu x +\nu y &\beta y  \\
\end{array}
$$
Note that $yx$ must be non-zero since the algebra is non-commutative.

Using the fact that $\mathfrak A$ is associative we derive the following
equations:
$$
\begin{array}{ccc}
x(yx)=(xy)x=0 &  \imply & \gamma+\alpha\mu=0 \\
(yx)y=y(xy)=0 &  \imply& \gamma+\beta\nu=0 \\
(yx)x=y(x^2) & \imply& \nu =\alpha \\
y(yx)=(y^2)x &  \imply& \mu=\beta\\
\end{array}
$$

Thus the table reduces to:
$$
\begin{array}{ccccc}
\times_{\mathfrak A}&\vline  & 1 &x &y  \\
\hline
1&\vline &1 &x &y  \\
x&\vline &x &\alpha x &0  \\
y&\vline &y &-\alpha\beta+\beta x +\alpha y &\beta y  \\
\end{array}
$$
We see that we only have three principally different cases - one when both
$\alpha$ and $\beta$ are non-zero ($2\times 2$ upper-triangular matrices):
$$
\begin{array}{ccccc}
\times_{\mathfrak A}&\vline  & 1 &x &y  \\
\hline
1&\vline &1 &x &y  \\
x&\vline &x &x &0  \\
y&\vline &y &-1+ x + y &y  \\
\end{array}
$$
and another two when one of $\alpha$ and $\beta$ vanishes:
$$
\begin{array}{ccccc}
\times_{\mathfrak A}&\vline  & 1 &x &y  \\
\hline
1&\vline &1 &x &y  \\
x&\vline &x &0 &0  \\
y&\vline &y & x & y  \\
\end{array}
$$
and
$$
\begin{array}{ccccc}
\times_{\mathfrak A}&\vline  & 1 &x &y  \\
\hline
1&\vline &1 &x &y  \\
x&\vline &x & x &0  \\
y&\vline &y & y &0 \\
\end{array}
$$
Both of the above two cases give valid multiplication tables which are, in fact,
isomorphic to the first case ($2\times2$ upper-triangular matrices), however,
such tables are not possible with our choice of $x$, $y$ and $F$. (Indeed, 
for these tables characteristic polynomial evaluated in $F$ would give zero -
which contradicts our assumption that $F$ is generic).

Hence there is only one non-commutative $3$-dimensional associative algebra with
unity.
\section{Jordan theory in multiplication table analysis}

While the methods presented above worked fine for algebras of small dimensions,
one cannot expect to find the direct sum of $\Stab_F(\alpha)$ equal to the algebra
$\mathfrak A$ in the general case. In the following we will assume that characteristic
polynomial of the algebra $\mathfrak A$ is non-zero. Fix a generic element $F$ of
$\mathfrak A^*$. Let $A$ be the multiplication table of $\mathfrak A$
evaluated in $F$. Let $\mu$ be such a value
that $A-\mu A^T$ is invertible.

\begin{definition}
We define $V_{-1}(\alpha)=\left\{0\right\}$.
\end{definition}

\begin{definition}
We introduce $V_0(\alpha)=\Stab_F(\alpha)$.
\end{definition}

\begin{definition}\label{V_k}
We define $V_k(\alpha)$ - a space of "Jordan vectors" - as
$$
V_{k+1}(\alpha):=\left\{b\in\mathfrak A:\exists a\in V_k(\alpha) \imply\left(\forall 
		x\in\mathfrak A\imply F(bx)-\alpha F(xb)=F(ax)-\mu F(xa)\right)\right\}
$$
or in terms of multiplication table of $\mathfrak A$:
$$
V_{k+1}(\alpha):=\left\{b\in\mathfrak A:\exists a\in V_k(\alpha) \imply
    bA-\alpha bA^T= a A -\mu a A^T \right\}
$$
\end{definition}

\begin{statement} The space $V_k(\alpha)$ does not depend on the choice of $\mu\neq\alpha$.
\end{statement}
\begin{proof}
Indeed, for $a\in V_k(\alpha)$ and $\mu_1$ and $\mu_2$ different from
$\alpha$:
\begin{eqnarray*}
\lefteqn{F(ax)-\mu_2 F(xa)=}\cr
&&\quad=(1-\delta)F(ax)+\delta\alpha F(xa) + \delta (F(a'x)-\mu_1F(xa'))
    -\mu_2 F(xa)=\cr
&&\quad=(1-\delta)F(ax)+(\delta\alpha-\mu_2)F(xa)+\delta(F(a'x)-\mu_1F(xa'))=\cr
&&\quad=F(((1-\delta)a+\delta a')x)-\mu_1 F(x ( (\delta\alpha-\mu_2)a/\mu_1 -\delta a'))
\end{eqnarray*}
Choosing $\delta=\frac{\mu_1+\mu_2}{\mu_1+\alpha}$ and $a''=(1-\delta)a+\delta a'$
We get $F(ax)-\mu_2 F(xa)=F(a''x)-\mu_1 F(xa'')$, where $a''$ is also in $V_{k-1}(\alpha)$.
A similar argument is used to prove the statement when one of $\mu_1$, $\mu_2$ is
infinity.
\end{proof}

 For $V_k(\infty)$ we need a special definition:
\begin{definition}
$$
V_{k+1}(\infty):=\left\{b\in\mathfrak A:\exists a\in V_k(\alpha) \imply\left(\forall 
		x\in\mathfrak A\imply F(xb)=F(xa)-\mu F(xa)\right)\right\}
$$
or
$$
V_{k+1}(\infty):=\left\{b\in\mathfrak A:\exists a\in V_k(\alpha) \imply
    bA^T= a A -\mu a A^T \right\}
$$
\end{definition}

\begin{statement}\label{V_k_prop} The spaces $V_k(\alpha)$ possess the following properties:
\begin{itemize}
\item $V_k(\alpha) \subset V_{k+1}(\alpha)$
\item $V_k(\alpha) \cdot V_m(\beta) \subset V_{k+m}(\alpha\beta)$
\end{itemize}
($\alpha \neq \infty$, $\beta \neq \infty$)
\end{statement}
\begin{proof}
The first part of the statement follows by induction from the fact that
$V_{-1}(\alpha)\subset V_0(\alpha)$.

To prove the second part we will make induction on the parameter $N=k+m$.
The base of induction follows immediately from properties of $\Stab_F(\alpha)$ 
(statement \ref{stab_times}).
Assume that the statement is true for all $k$ and $m$ such that $k+m<N+1$.
Let us pick certain $\alpha$ and $\beta$. Let $b_1\in V_k$, $b_2 \in V_m$, where
both $k+m=N+1$. Let $a_1\in V_{k-1}(\alpha)$ be an element corresponding to $b_1$
according to definition \ref{V_k} and $a_2\in V_{m-1}(\beta)$ be the element corresponding to $b_2$.
Let $x$ be an arbitrary element of $\mathfrak A$. Assume $\mu=\infty$ (since 
$\alpha$ and $\beta$ are finite).
Then:
\begin{eqnarray*}
\lefteqn{F(b_1b_2x)=}&&\cr
&&\qquad=\alpha F(b_2xb_1) + F(b_2xa_1)=\cr
&&\qquad=\alpha\beta F(xb_1b_2)+\alpha F(xb_1a_2)+F(b_2xa_1)=\cr
&&\qquad=\alpha\beta F(xb_1b_2)+F(x(\alpha b_1 a_2+\beta a_1 b_2+a_1a_2))
\end{eqnarray*}
Now by assumption of induction we have 
$$
\alpha b_1 a_2+\beta a_1 b_2+a_1a_2 \subset V_{k+m-1}(\alpha\beta)
$$
which concludes the proof.
\end{proof}

\begin{statement} For all $\alpha \neq 0$ we have
\begin{eqnarray*}
V_k(\alpha)\cdot V_m(\infty)&\subset& V_{k+m}(\infty)\cr
V_k(\infty)\cdot V_m(\alpha)&\subset& V_{k+m}(\infty)
\end{eqnarray*}
\end{statement}

\begin{statement} 
$$
V_k(0)\cdot V_m(\infty)\subset V_k(0)\intersect V_m(\infty)
$$
\end{statement}

The proof of the above two statements is analogous to proof of statement \ref{V_k_prop}.

\begin{statement}\label{kernel F}Let $\mathfrak A$ be an associative algebra with unity. Then
$F(V_k(\alpha))=0$ for all $\alpha \neq 1$.
\end{statement}
\begin{proof}
Case $k=0$: for all $b\in V_0(\alpha)$ we must have
$$
F(bx)=\alpha F(xb)
$$
Setting $x=1$ we get $F(b)=\alpha F(b)$, hence $F(b)=0$.

For arbitrary $k$ we proceed by induction. Again let us set $x=1$ in the definition
\ref{V_k}. We get
$$
F(b)-\alpha F(b)=F(a)
$$
But we already know that $F(a)=0$, thus $F(b)=0$ as well.

The case $\alpha=\infty$ is similar:  from the definition we get
$$
F(b)=F(a)
$$
and from the definition of $V_0(\infty)$ it follows that $F(V_0(\infty))=0$.
\end{proof}

%
%
%
%

\begin{statement} \ 
\begin{itemize}
\item $V_k(\alpha)\subset V_{k+1}(\alpha)$
\item $V_k(\alpha)$ stabilizes for sufficiently large $k$.
\end{itemize}
\end{statement}

\begin{statement}
$$
\mathfrak A=\oplus_{\alpha} V_N(\alpha)
$$
(here N is sufficiently large)
\end{statement}
\begin{proof}
We will reduce this problem to the classical theorem about Jordan normal form 
for matrices.

Recall the definition of $V_{k+1}(\alpha)$:
$$
V_{k+1}(\alpha):=\left\{b\in\mathfrak A:\exists a\in V_k(\alpha) \imply
    bA-\alpha bA^T= a A -\mu a A^T \right\}
$$

The equation $bA-\alpha bA^T= a A -\mu a A^T$ can be transformed as follows:
\begin{eqnarray*}
\lefteqn{bA-\alpha bA^T= a A -\mu a A^T \equivalent}\cr
&&\qquad \equivalent b\left( \left[A-\mu A^T\right]-(\alpha-\mu)A^T\right)=a(A-\mu A^T)\equivalent\cr
&&b \left(1-(\alpha-\mu)A^T(A-\mu A^T)^{-1}\right)=a\equivalent
\end{eqnarray*}

Thus $V_k(\alpha)=\ker \left(1-(\alpha-\mu)A^T(A-\mu A^T)^{-1}\right)^k$. 
Hence $V_k(\alpha)$ 
is nothing more but the space of Jordan vectors of rank $k$ corresponding to 
eigenvalue $(\alpha-\mu)^{-1}$ (note also a non-standard order of multiplication
between vectors and operators).
\end{proof}

Let us summarize our findings in the following theorem:
\begin{theorem}\label{V_n} Let $\mathfrak A$ be an associative algebra over an algebraicly
closed field. Assume that characteristic polynomial of $\mathfrak A$ is non-zero.
Let $F$ be a generic element of $\mathfrak A^*$.
Then
\begin{itemize}
\item $\mathfrak A=\oplus V_N(\alpha)$, ($N=\dim \mathfrak A$)
\item $\dim V_N(\alpha)=\dim V_N(1/\alpha)$
\item $F(V_k(\alpha)\cdot V_m(\beta))=0$ when $\alpha \beta\neq 1$ (provided $1\in \mathfrak A$)
\item $F(V_k(\alpha))=0$ when $\alpha \neq 1$ (provided $1\in \mathfrak A$)
\item the multiplication relations for $V_N(\alpha)$ can be summarized in
       the following table:
$$
\begin{array}{ccccc}
& \vline & V_k(0) & V_k(\alpha) & V_k(\infty) \\
\hline
V_m(0) & \vline & V_{k+m}(0)  & V_{k+m}(0)  & V_{m}(0) \intersect V_{k}(\infty)=0   \\
V_m(\alpha) & \vline & V_{k+m}(0)  & V_{k+m}(\alpha\beta)  & V_{k+m}(\infty)   \\
V_m(\infty) & \vline & *  & V_{k+m}(\infty)  &  V_{k+m}(\infty)   
\end{array}
$$
\item $V_N(0)$, $V_N(1)$ and $V_N(\infty)$ are closed under multiplication
\item $V_N(0)$ is solvable (as Lie algebra)
\item $(\lambda-\alpha\gamma)^{\dim V_k(\alpha)}$ is the highest power of
$\lambda -\alpha \gamma$ that divides characteristic polynomial of $\mathfrak A$
(for $\alpha=\infty$ we use $\gamma$ in instead of $\lambda -\alpha \gamma$).
\end{itemize}
\end{theorem}
\begin{proof}
Most of these statements have already been discussed. The fact that
$V_m(0)\intersect V_k(\infty)=0$ follows by induction from $\Nil_F=0$ - 
which, in turn, is due to the requirement that characteristic polynomial
is non-zero. $\dim V_N(\alpha)=\dim V_N(1/\alpha)$ in view of the symmetry
$A \leftrightarrow A^T$, $\alpha \leftrightarrow 1/\alpha$. Solvability
of $V_N(1)$ is due to commutativity of $\Stab_F(1)$. The last statement
follows from the following computation: 
\begin{eqnarray*}
\lefteqn{\det\left(\lambda A +\gamma A^T\right)=}\\
&&\qquad=\det\left(\lambda + (\gamma+\mu\lambda)A^T(A-\mu A^T)^{-1}\right) \cdot \det(A-\mu A^T)
\end{eqnarray*}
\end{proof}

\begin{statement}\label{A_k}Let $\alpha\notin \left\{0,1,\infty\right\}$.
 Define $A_{\alpha}:V_k(\alpha)\times V_m(1/\alpha)\rightarrow \mathbb C$ as
$$
A_{\alpha}(a,b)=F(ab)
$$
then the map $A_{\alpha}^T$:
$$
A_{\alpha}^T(a,b)=F(ba)
$$
is related to $A_k$ by the following equation:
$$
A_{\alpha}(a,b)=A_{\alpha}^T(U_ka,b)
$$
where $U_k$  is an automorphism of $V_k(\alpha)$ with only eigenvalue $\alpha$.
\end{statement}
\begin{proof}
We will proceed by induction. For $k=0$ the statement is true with $U_0=\alpha$.
Assume that the statement holds for all $k<K$. By definition for $b\in V_{K+1}(\alpha)$
and $x\in V_{m}(1/\alpha)$ we have
$$
F(bx)=\alpha F(xb)+ F(xa)=F(x(\alpha b+a))
$$
where $a\in V_{K}(\alpha)$. Let us choose a basis $\left\{e_i\right\}$ in $V_{K+1}(\alpha)$
such that first $l=\dim V_K(\alpha)$ vectors are in $V_K(\alpha)$. Let $\left\{a_i\right\}$
be vectors $a_i \in V_K(\alpha)$ corresponding to vectors $e_{l+1}...e_{n}$. 
Then the operator $U_{K+1}$ is given by the matrix
$$
\begin{array}{cccc}
& \vline & e_1...e_l & e_{l+1}...e_n \\
\hline
e_1...e_l & \vline & U_{K}& \left(a_{l+1}...a_{n}\right)\\
e_{l+1}...e_n & \vline &0 &\alpha \cdot 1
\end{array}
$$
\end{proof}

\begin{statement}\label{V_n duality} Let $\mathfrak A$ be an associative algebra with non-zero 
characteristic polynomial.
Then the map $A_{\alpha}$ defined in statement \ref{A_k} is a duality between spaces
$V_N(\alpha)$ and $V_N(1/\alpha)$.
\end{statement}
\begin{proof}
Indeed, consider the multiplication table evaluated in point $F$:
$$
\begin{array}{ccccccc}
 & \vline & V_N(0) & V_N(\alpha) & V_N(1) & V_N(1/\alpha) & V_N(\infty) \\
 \hline
V_N(0) & \vline  &   0  &  0   &  0   &   0   & 0\\
V_N(\alpha) & \vline  &  0   &  0   &  0   &  U(\alpha)A_{\alpha}^T    & 0\\
V_N(1) & \vline  &  0   &  0   & Q_F    &  0    &0 \\
V_N(1/\alpha) & \vline  &  0   &  A_{\alpha}^T   &  0  &   0   & 0\\
V_N(\infty) & \vline  &  A_{\infty}   &  0   &  0   &   0   & 0\\
\end{array}
$$
(here $A_{\infty}(a,b)=F(ab)$ is a bilinear form defined on $V_N(\infty)\times V_N(0)$)
The characteristic polynomial is then given by the following formula:
\begin{eqnarray*}
\lefteqn{\chi(\lambda,\mu,F)=\det(\lambda A_{\infty})\det(\mu A_{\infty})\det((\lambda+\mu)Q_F)\cdot}\\
&&\qquad\cdot\prod_{\alpha_k\notin \left\{0,1,\infty\right\}}{\det((\lambda U(\alpha_k)+\mu)A_{\alpha_k}^T)}
\end{eqnarray*}
Thus the only way for $\chi(\lambda,\mu,F)$ to be nonzero is for $A_{\infty}$, $Q_F$
and $A_{\alpha_k}$ to be invertible. Note, however, that the form $Q_F$ here is defined
on $V_N(1)$ - as opposed to $\Stab_F(1)$ as used in statement \ref{Q_F}.						 
\end{proof}

A consequence of the above results and classification of $3$-dimensional algebras
is the following theorem:
\begin{theorem}\label{V_alpha is zero} Let $\mathfrak A$ be an associative algebra with unity and
Lie index $1$. Then $\dim V_k(\alpha)=0$ for all $\alpha$ except $0$, $1$ and 
$\infty$.
\end{theorem}
\begin{proof}
Assume the converse, i.e. there exists $\alpha\notin \left\{0,1,\infty\right\}$
such that 
$$
\dim V_N(\alpha)>0
$$
If there is more than one candidate let us choose one that has maximum absolute
value. By theorem \ref{V_n} $\alpha^{-1}$ should have the smallest possible value
 (because otherwise there would have existed non-zero $V_N(1/\alpha_1)$ and due
 to symmetry $V_N(\alpha_1)$ as well). 
 
 Pick two vectors $v\in V_N(\alpha)$ and
 $u \in V_N(1/\alpha)$. We will show that their product is zero. Indeed, 
 $\mathbb C\cdot 1+\mathbb C \cdot v + \mathbb C \cdot u$ forms a $3$-dimensional
 subalgebra\footnote{
 Note that we are being slightly inaccurate. If it happens
that $\alpha=-1$ then vector spaces $V_N(\alpha)$ and $V_N(1/\alpha)$ coincide.
In this case it is possible for vectors $u$ and $v$ to be parallel. However,
it is straightforward to show that the product $u\cdot v$ must still be zero.
 } with unity in $\mathfrak A$. Due to classification of $3$-dimensional
 associative algebras with unity we know that they have trivial spaces $V_N(\alpha)$
 for $\alpha \notin \left\{0,1,\infty\right\}$. Thus 
 $v \cdot u \in V_N(1)=\mathbb C \cdot 1$ must
 be $0$. 
 
 This proves that the map $V_N(\alpha) \times V_N(1/\alpha)\rightarrow \mathbb C \cdot 1$ 
is trivial. However this contradicts statement \ref{V_n duality} which states that
characteristic polynomial must be zero in this case (and it is not in view of
statement \ref{Q_F}).
\end{proof}

\begin{definition} We define $\phi_{\alpha,\beta}: V_N(\alpha)\times V_N(\beta)
\rightarrow V_N(\alpha\beta):(x,y)\rightarrow xy$.
\end{definition}

\begin{theorem} The image of $\phi_{\alpha,1/\alpha}$ is an ideal in $V_N(1)$ for all
$\alpha\notin\left\{0,\infty\right\}$.
\end{theorem}
\begin{proof}
Indeed, $V_N(1)$ acts on $V_N(\alpha)$ for all $\alpha \notin\left\{0,\infty\right\}$
both from the left and right. The map $\phi_{\alpha,\beta}$ is eqivariant with respect
to left action on the first argument and right action on the second argument.
Consequently, the image of $\phi_{\alpha,\beta}$ is stable under both left and
right action of $V_N(1)$.
\end{proof}

\section{Associative algebras with unity and index $1$}
Due to theorem \ref{V_alpha is zero} we know that for an associative algebra
with unity and index $1$ the multiplication table is
$$
\begin{array}{ccccc}
 & \vline & x_2 & 1  & y_2 \\
 \hline
x_1 & \vline & x_1x_2  & x_1   &  0  \\
1 & \vline & x_2   &  1   &  y_2   \\
y_1  & \vline & -1\cdot A(y_1,x_2)+B(y_1,x_2)+C(y_1,x_2)    &  y_1   &  y_1y_2  
\end{array}
$$
Here $x_1, x_2, B(y_1,x_2) \in V_N(0)$ and $y_1, y_2, C(y_1,x_2) \in V_N(\infty)$.
Very much like in the case of $3$-dimensional algebras with unity we have equations
on $A$, $B$ and $C$ due to associativity:
\begin{eqnarray}
\label{a_b}x_1(yx_2)=(x_1y)x_2 & \implies &  -x_1 A(y,x_2)+x_1 B(y,x_2)=0\\
\label{a_c}(y_1x)y_2=y_1(xy_2) & \implies &  -y_2 A(y_1,x)+C(y_1,x)y_2=0     
\end{eqnarray}
$$
\begin{array}{ccl}
\lefteqn{(yx_1)x_2=y(x_1x_2)\implies}\\
&&\qquad \implies -x_2 A(y,x_1)+B(y_1,x_1)x_2-\\
&&\qquad\qquad -1\cdot A(C(y,x_1),x_2)+B(C(y,x_1),x_2)+C(C(y,x_1),x_2)=\\
&&\qquad=-1\cdot A(y,x_1x_2)+B(y,x_1x_2)+C(y,x_1,x_2)     
\end{array}
$$
$$
\begin{array}{ccl}
\lefteqn{y_1(y_2x)=(y_1y_2)x \implies }\\
&& \qquad \implies -y_1 A(y_2,x)+y_1 C(y_2,x)-\\
&&\qquad\qquad - 1\cdot A(y_1,B(y_2,x))+B(y_1,B(y_2,x))+C(y_1,B(y_2,x))=\\
&& \qquad = -1\cdot A(y_1y_2,x)+B(y_1y_2,x)+C(y_1y_2,x)
\end{array}
$$
An immediate consequence of these equations are the following properties of
$A$, $B$ and $C$:
\begin{eqnarray}
\label{a_x}A(y,x_1x_2)&=&A(C(y,x_1),x_2) \\
\label{b_x}B(y,x_1x_2)&=&B(C(y,x_1),x_2)+[B(y,x_1),x_2]\\
\label{c_x}C(y,x_1x_2)&=&C(C(y,x_1),x_2)\\
\label{a_y}A(y_1y_2,x)&=&A(y_1,B(y_2,x))\\
\label{b_y}B(y_1y_2,x)&=&B(y_1,B(y_2,x)) \\
\label{c_y}C(y_1y_2,x)&=&C(y_1,B(y_2,x))+[y_1,C(y_2,x)]
\end{eqnarray}
Since, by statement \ref{kernel F}, $F(V_N(\alpha))=0$ when $\alpha\neq1$ we must have that
$A(y,x)$ is non-degenerate - and thus establishes duality between $V_N(0)$ and 
$V_N(\infty)$. Property \ref{c_x} states that $C(y,x)$ establishes a right action
of $V_N(0)$ on $V_N(\infty)$ and property \ref{a_x} states that this action is
the same one as induced from multiplication on $V_N(0)$ 
by duality $A$. Analogously, property \ref{b_y} states
that $B$ describes  left action of $V_N(\infty)$ on $V_N(0)$ and property \ref{a_y}
states that this is the same action as induced from $V_N(\infty)$ by duality $A$. Thus
$$
\mathfrak A=\mathfrak H \oplus \mathbb C\cdot 1 \oplus \mathfrak H'
$$
where both $\mathfrak H:=V_N(0)$ and $\mathfrak H'$  possess operations of 
associative multiplication. The operation of multiplication on $\mathfrak H'$
gives rise to a map $\Delta:\mathfrak H\rightarrow \mathfrak H \otimes \mathfrak H$:
$$
\langle\Delta(x),y_1\otimes y_2\rangle=\langle x, y_1y_2\rangle
$$
Equation \ref{a_b} states:
$$
- x_1 A(y_1,x_2)+x_1 B(y_1,x_2)=0
$$
Or, equivalently (using the fact that $A$ is a duality),
$$
A(y_2,x_1) A(y_1,x_2)=A(y_2,x_1 B(y_1,x_2))
$$
Rewriting in terms of $\mathfrak H$ and $\mathfrak H'$:
\begin{eqnarray*}
\lefteqn{\langle x_1,y_2\rangle \langle x_2, y_1\rangle=\langle x_1B(y_1,x_2),y_2\rangle\equivalent}\\
&&\qquad\equivalent \langle x_1,y_2\rangle \langle x_2, y_1\rangle=\langle B(y_1,x_2), C(y_2,x_1)\rangle
\end{eqnarray*}
Equation \ref{a_c} states:
$$
- y_2 A(y_1,x_2)+C(y_1,x_2)y_2=0   
$$
or, equivalently (using the fact that $A$ is a duality),
$$
A(y_2,x_1)A(y_1,x_2)=A(C(y_1,x_2)y_2,x_1)
$$
Rewriting in terms of $\mathfrak H$ and $\mathfrak H'$ we obtain:
\begin{eqnarray*}
\lefteqn{\langle x_1, y_2\rangle\langle x_2,y_1\rangle=\langle x_1,C(y_1,x_2)y_2\rangle\equivalent}\\
&&\qquad \equivalent \langle x_1, y_2\rangle\langle x_2,y_1\rangle=\langle B(y_2,x_1), C(y_1,x_2)\rangle
\end{eqnarray*}
Hence equations \ref{a_b} and \ref{a_c} follow from identity
\begin{equation}\label{<xy><xy>}
\langle x_1,y_2\rangle \langle x_2, y_1\rangle=\langle B(y_1,x_2), C(y_2,x_1)\rangle
\end{equation}
and properties \ref{a_x}, \ref{c_x}, \ref{a_y} and \ref{b_y}.

Equation \ref{b_x} states:
$$
B(y_2,x_1x_2)=B(C(y_2,x_1),x_2)+[B(y_2,x_1),x_2]
$$
or, equivalently (using the fact that $A$ is a duality),
$$
\langle B(y_2,x_1x_2),y_1\rangle=\langle B(C(y_2,x_1),x_2),y_1\rangle +
\langle B(y_2,x_1)x_2,y_1\rangle - \langle x_2 B(y_2,x_1),y_1\rangle
$$
After transformation we get
$$
\langle x_1x_2,y_1y_2\rangle=\langle x_2,y_1C(y_2,x_1)\rangle +
\langle B(y_2,x_1)x_2,y_1\rangle - \langle B(y_2,x_1),C(y_1,x_2)\rangle
$$

Equation \ref{c_y} states:
$$
C(y_1y_2,x_1)=C(y_1,B(y_2,x_1))+[y_1,C(y_2,x_1)]
$$
or, equivalently (using the fact that $A$ is a duality),
$$
\langle x_2, C(y_1y_2,x_1)\rangle=\langle x_2, C(y_1,B(y_2,x_1))\rangle +
\langle x_2,y_1 C(y_2,x_1) \rangle - \langle x_2,C(y_2,x_1)y_1 \rangle
$$
After transformation we get
$$
\langle x_1x_2,y_1y_2\rangle=\langle B(y_2,x_1)x_2,y_1\rangle +
\langle x_2,y_1 C(y_2,x_1) \rangle - \langle B(y_1,x_2),C(y_2,x_1) \rangle
$$
which is the same  one as for equation \ref{b_x}. 

Thus we  have reduced the original equations imposed by associativity to the
following conditions:
\begin{itemize}
\item $V_N(0)=\mathfrak H$ is an associative algebra with additional operation
of multiplication on its dual space. $V_N(\infty)=\mathfrak H'$ is dual to $V_N(0)$
\item $A=\langle x,y\rangle$ is given by duality between $\mathfrak H$ and $\mathfrak H'$
\item $B$ is given by the dual of the right action by multiplication on $\mathfrak H'$
$$
\langle B(y_2,x),y_1\rangle=\langle x,y_1y_2\rangle
$$
\item $C$ is given by the dual of the left action by multiplication on $\mathfrak H$
$$
\langle x_2,C(y,x_1)\rangle=\langle x_1x_2,y\rangle
$$
\item there are two additional conditions on multiplication structures on $\mathfrak H$
and $\mathfrak H'$:
\begin{equation}\label{rank1_eqn}
\langle x_1,y_1\rangle \langle x_2, y_2\rangle=\langle B(y_1,x_1), C(y_2,x_2)\rangle
\end{equation}
and
\begin{equation}\label{homo_eqn}
\langle x_1x_2,y_1y_2\rangle=\langle x_2,y_1C(y_2,x_1)\rangle +
\langle B(y_2,x_1)x_2,y_1\rangle - \langle B(y_1,x_2),C(y_2,x_1)\rangle
\end{equation}
\end{itemize}

\begin{definition} We define $B_R=B$, $C_L=C$ and
$$
\langle B_L(y_1,x),y_2\rangle=\langle x , y_1y_2\rangle
$$
$$
\langle x_1,C_R(y,x_2)\rangle=\langle x_1x_2,y\rangle
$$
\end{definition}

In this notation and using equation \ref{rank1_eqn} one can rewrite equation \ref{homo_eqn} 
as
\begin{eqnarray}\label{homo_eqn.2}
\lefteqn{\langle x_1x_2,y_1y_2\rangle+\langle x_1,y_1\rangle \langle x_2, y_2\rangle=}\nonumber\\
&\qquad=\langle B_L(y_1,x_2),C_L(y_2,x_1)\rangle +
\langle B_R(y_2,x_1) ,C_R(y_1,x_2)\rangle
\end{eqnarray}

In this form this equation is clearly symmetric with respect to transformation $x \leftrightarrow y$.

\if \extended

We will now study equation \ref{rank1_eqn} in more detail. 
Let $U_0=\image B_R=\image B$ and $U_1=\image C_L=\image C$. Let $W_0=U_0 \intersect U_1^{\orthogonal}$
and $W_1=U_1 \intersect U_0^{\orthogonal}$.
Since $x_1$, $y_1$, $x_2$ and
$y_2$ are arbitrary we must have $\dim U_0/W_0=1$ and
$\dim U_1/W_1=1$. Let us choose elements $b_0\in U_0$ and $c_0 \in U_1$ such
that
$$
B(y,x)=b_0 \langle y,x \rangle+B_1(y,x)
$$
and
$$
C(y,x)=c_0\langle y,x\rangle+C_1(y,x)
$$
where $B_1(y,x)\in W_0$ and $C_1(y,x)\in W_1$. Substituting these expressions
into equation \ref{rank1_eqn} yields $\langle c_0,b_0\rangle=1$, $\langle c_0,W_0\rangle=0$ and
$\langle W_1,b_0\rangle=0$.

Since $B$ is a dual of right action by $y$ we get
$$
\langle x,y_1y_2\rangle=\langle B(y_2,x),y_1\rangle=\langle b_0 \langle y_2,x\rangle+
B_1(y_2,x),y_1\rangle
$$
Using $y_1=c_0$ we obtain
$$
\langle x,c_0 y_2\rangle = \langle x, y_2\rangle
$$
and thus $c_0y=y$. Analogously, one proves that $xb_0=x$.

\fi

\references

\end{document}